\documentclass[12pt]{article}

\usepackage{epsf,psfrag,amsmath,amsfonts, amssymb, bbm, harvard}

\begin{document}
\bibliographystyle{harvard}

\title{Central Limit Theorem and convergence to stable laws in Mallows 
distance}
\author{Oliver Johnson and Richard Samworth \\ 
Statistical Laboratory, 
University of Cambridge, \\
Wilberforce Road, Cambridge, CB3 0WB, UK. }
\date{\today}
\maketitle

Running Title: CLT and stable convergence in Mallows distance


\newtheorem{theorem}{Theorem}[section]
\newtheorem{lemma}[theorem]{Lemma}
\newtheorem{proposition}[theorem]{Proposition}
\newtheorem{corollary}[theorem]{Corollary}
\newtheorem{conjecture}[theorem]{Conjecture}
\newtheorem{definition}[theorem]{Definition}
\newtheorem{example}[theorem]{Example}
\newtheorem{remark}[theorem]{Remark}
\newtheorem{condition}{Condition}
\newtheorem{main}{Theorem}
\newtheorem{assumption}[theorem]{Assumption}
\setlength{\parskip}{\parsep}

\def \outlineby #1#2#3{\vbox{\hrule\hbox{\vrule\kern #1%
\vbox{\kern #2 #3\kern #2}\kern #1\vrule}\hrule}}%
\def \endbox {\outlineby{4pt}{4pt}{}}%
\newenvironment{proof}
{\noindent{\bf Proof\ }}{{\hfill \endbox
}\par\vskip2\parsep}
\newenvironment{pfof}[2]{\removelastskip\vspace{6pt}\noindent
 {\it Proof  #1.}~\rm#2}{\par\vspace{6pt}}

\hfuzz20pt

\newcommand{\Section}[1]{\setcounter{equation}{0} \section{#1}}
\newcommand{\var}{{\rm{Var\,}}}
\newcommand{\cov}{{\rm{Cov\,}}}
\newcommand{\tends}{\rightarrow \infty}
\newcommand{\C}{{\mathcal{C}}}
\newcommand{\ep}{{\mathbb {E}}}
\newcommand{\pr}{{\mathbb {P}}}
\newcommand{\sgn}{{\rm{sgn}}}
\newcommand{\re}{{\mathbb {R}}}
\newcommand{\I}{\mathbbm{1}}
\newcommand{\vc}[1]{{\bf {#1}}}
\newcommand{\conpr}{\buildrel{\pr}\over\longrightarrow}
\newcommand{\cond}{\stackrel{d}{\rightarrow}}
\newcommand{\condo}{\stackrel{d^{\circ}}{\rightarrow}} 
\newcommand{\q}{{\cal{Q}}}
\newcommand{\blah}[1]{}
\begin{abstract}
We give a new proof of the classical Central Limit Theorem, in
the Mallows ($L^r$-Wasserstein) distance. Our proof is elementary in the
sense that it does not require complex analysis, but rather makes use of
a simple subadditive inequality related to this metric. The key is to 
analyse the case where equality holds. We provide some results concerning
rates of convergence. We also consider convergence to stable distributions,
and obtain a bound on the rate of such convergence.
\end{abstract}
\section{Introduction and main results} \label{sec:intro}
The spirit of the Central  Limit Theorem, that  normalised sums of
independent random variables converge to a normal distribution, can be
understood  in   different  senses,  according   to  the distance
used. 
For example, in addition to the standard Central Limit Theorem 
in the sense of weak convergence, we
mention the proofs in \citeasnoun{prohorov} of $L^1$ convergence of densities,
in \citeasnoun{gnedenko} of $L^{\infty}$ convergence
of densities, in \citeasnoun{barron} of convergence in 
relative entropy and in \citeasnoun{shimizu} and
\citeasnoun{johnson5} of convergence in Fisher information.

In this  paper we consider the Central
Limit Theorem with respect to the Mallows distance and 
prove convergence to stable laws in the infinite
variance setting. We study the rates of convergence in both cases.
\begin{definition} \label{def:mallows}
For any $r > 0$, we define the
Mallows  $r$-distance between  
probability distribution functions  $F_X$ and $F_Y$  as
$$ d_r(F_X,F_Y)  = \left(  \inf_{(X,Y)} \ep |  X - Y  |^r \right)^{1/r},$$
where the  infimum is  taken over pairs $(X,Y)$ 
whose marginal distribution functions are $F_X$ and
$F_Y$ respectively, and may be infinite.
Where it causes no confusion, we write $d_r(X,Y)$ for $d_r(F_X,F_Y)$.
\end{definition}
Define 
$\mathcal{F}_r$ to be the set of distribution functions $F$ such
that $\int |x|^r dF(x) < \infty$. \citeasnoun{bickel} show that 
for $r \geq 1$, $d_r$ is a metric on $\mathcal{F}_r$.
If $r < 1$, then $d_r^r$ is a metric on
$\mathcal{F}_r$. In considering stable
convergence, we shall also be concerned with the case where the absolute
$r$th moments are not finite.

Throughout the paper, we write $Z_{\mu,\sigma^2}$ for a $N(\mu,\sigma^2)$
random variable, $Z_{\sigma^2}$ for a $N(0,\sigma^2)$ random variable,
and $\Phi_{\mu,\sigma^2}$ and $\Phi_{\sigma^2}$ for their respective
distribution functions.
We establish the following main theorems:
\begin{theorem} \label{thm:main}
Let $X_1, X_2, \ldots $ be independent and identically distributed random 
variables with mean zero and finite variance $\sigma^2 > 0$, and let 
$S_n = (X_1 + \ldots + X_n)/\sqrt{n}$.  Then
$$ \lim_{n \tends} d_2(S_n, Z_{\sigma^2}) = 0.$$
\end{theorem}
Moreover, Theorem \ref{thm:mainr} shows that for any $r \geq 2$, 
if $d_r(X_i, Z_{\sigma^2}) < \infty$,  
then $ \lim_{n \tends} d_r(S_n, Z_{\sigma^2}) = 0.$
Theorem \ref{thm:main} 
implies the standard Central Limit Theorem in the sense
of weak convergence \cite[Lemma 8.3]{bickel}.
\begin{theorem} \label{thm:mainstab} Fix $\alpha \in (0,2)$,
and let $X_1, X_2, \ldots $ be independent random 
variables (where $\ep X_i =0$, if $\alpha > 1$), and  
$S_n = (X_1 + \ldots + X_n)/n^{1/\alpha}$.  If there exists an
$\alpha$-stable random variable $Y$ such that
$\sup_i d_{\beta}(X_i, Y) < \infty$ for some $\beta \in (\alpha,2]$,
then $\lim_{n \tends} d_{\beta}(S_n, Y) = 0$. In fact 
$$d_{\beta}(S_n, Y) \leq \frac{2^{1/\beta}}{n^{1/\alpha}} 
\left( \sum_{i=1}^n d^{\beta}_{\beta}(X_i,Y) \right)^{1/\beta},$$
so in the identically distributed  
case the rate of convergence is $O(n^{1/\beta - 1/\alpha})$.
\end{theorem}
See also Rachev and R{\"u}schendorf (1992,1994), who obtain
similar results using different techniques in the case of identically
distributed $X_i$  and strictly symmetric $Y$. 
In Lemma \ref{lem:dsna} we exhibit a large class ${\cal C}_K$ of
distribution functions $F_X$ 
for which $d_{\beta}(X,Y) \leq K$, 
so the theorem can be applied.

Theorem \ref{thm:main} follows by understanding the subadditivity of 
$d_2^2(S_n, Z_{\sigma^2})$ (see Equation (\ref{eq:sub})). We consider
the powers-of-two subsequence $T_k = S_{2^k}$, and
use R\'{e}nyi's method, introduced in
\citeasnoun{renyi2} to provide a proof of convergence to equilibrium of Markov 
chains; see also \citeasnoun{kendall}. This technique was also
used in \citeasnoun{csiszar2} to show 
convergence to Haar measure for convolutions of measures on compact groups,
and in \citeasnoun{shimizu} to show convergence of Fisher information
in the Central Limit 
Theorem. The method has four stages:
\begin{enumerate}
\item{Consider independent and identically distributed random variables $X_1$ 
and $X_2$ with mean $\mu$ and variance $\sigma^2 > 0$, and write $D(X)$ for 
$d_2^2(X, Z_{\mu,\sigma^2})$. 
In Proposition \ref{prop:iff}, we observe that 
\begin{equation} \label{eq:subadd}
D \left( \frac{ X_1 + X_2}{\sqrt{2}} \right) \leq D(X_1),
\end{equation}
with equality if and only if $X_1,X_2 \sim Z_{\mu,\sigma^2}$.
Hence $D(T_k)$ is decreasing and bounded below, so converges to some $D$.}
\item{In Proposition \ref{prop:rich}, we use a compactness argument
to show that there exists a strictly increasing sequence $k_r$ and a random 
variable $T$ such that
$$ \lim_{r \tends} D(T_{k_r}) = D(T).$$
Further,
$$\lim_{r \tends} D(T_{k_r+1}) = \lim_{r \tends} D \left( 
\frac{T_{k_r} + T_{k_r}'}{\sqrt{2}} \right) = 
D \left( \frac{T + T'}{\sqrt{2}} \right),$$
where the $T'_{k_r}$ and $T'$ are independent copies of $T_{k_r}$ and $T$ 
respectively.}
\item{We combine these two results: since $D(T_{k_r})$ and
$D(T_{k_r +1})$ are both subsequences of the convergent subsequence $D(T_k)$,
they must have a common limit. That is,
$$ D = D(T) = D \left( \frac{ T + T'}{\sqrt{2}} \right),$$
so by the condition for equality in Proposition \ref{prop:iff}, we 
deduce that $T \sim N(0,\sigma^2)$ and $D =0$.}
\item{Proposition \ref{prop:iff} implies the standard subadditive relation
$$ (m+n) D(S_{m+n}) \leq m D(S_m) + n D(S_n).$$
Now Theorem 6.6.1 of  \citeasnoun{hille}
implies that $D(S_n)$ converges to $\inf_n D(S_n) = 0$.}
\end{enumerate}
The proof of Theorem \ref{thm:mainstab} is given in Section \ref{sec:stable}.
\section{Subadditivity of Mallows distance}
The Mallows distance and related metrics originated
with a transportation problem posed by Monge in 1781  
\cite[pp.329--330]{rachev2,dudley}.  
Kantorovich generalised this problem, and
considered the distance obtained by minimising $\ep c(X,Y)$, for a
general metric $c$ (known as the cost function), over all joint
distributions of pairs 
$(X,Y)$ with fixed marginals.  This distance is also known
as the Wasserstein metric. \citeasnoun{rachev2} reviews
applications to differential geometry, 
infinite-dimensional linear programming and information theory, 
among many others.  
\citeasnoun{mallows} 
focused on the metric which we have called $d_2$, while $d_1$
is sometimes called the Gini index.

In Lemma \ref{lem:covvalue} below,  we review  the existence and 
uniqueness of the construction which
attains the infimum in Definition \ref{def:mallows}, using the concept
of a quasi-monotone function.
\begin{definition}
A function $k: \re^2 \rightarrow \re$ 
induces a signed measure $\mu_k$ on $\re^2$ given by
$$ \mu_k \left\{ (x, x'] \times (y,y'] \right\}
= k(x,y) + k(x',y') - k(x,y') - k(x',y).$$
We say that $k$ is quasi-monotone if $\mu_k$ is a non-negative measure.
\end{definition}
The function $k(x,y) = -|x-y|^r$ is quasi-monotone for $r \geq 1$,
and if $r > 1$ then the measure $\mu_k$ is absolutely
continuous, with a density
which is  positive Lebesgue almost everywhere.
\citeasnoun[Corollary 2.1]{tchen} gives  the following 
result, a two-dimensional version of integration by parts.
\begin{lemma} \label{lem:intprt}
Let $k(x,y)$ be a quasi-monotone function and
let $H_1(x,y)$ and $H_2(x,y)$ be distribution
functions with the same marginals, where $H_1(x,y) \leq H_2(x,y)$ 
for all $x,y$. Suppose there exists an $H_1$- and $H_2$- integrable 
function $g(x,y)$, bounded on compact sets, such that 
$k(x^B,y^B) \leq g(x,y)$, where $x^B = (-B) \vee x \wedge B$. 
Then
$$ \int k(x,y) dH_2(x,y) - \int k(x,y) dH_1(x,y) 
= \int \left\{ H^{-}_2(x,y) - H^{-}_1(x,y) \right\} d\mu_k(x,y).$$
Here $H^{-}_i(x,y) = \pr( X < x, Y < y)$, where $(X,Y)$ have joint
distribution function $H_i$.
\end{lemma}
\begin{lemma} \label{lem:covvalue}
For $r \geq 1$, consider the joint 
distribution of pairs $(X,Y)$ where $X$ and $Y$ have fixed marginals 
$F_X$ and $F_Y$, both in $\mathcal{F}_r$. Then 
\begin{equation} \label{eq:covmax} \ep |X-Y|^r \geq \ep |X^* - Y^*|^r,
\end{equation}
where $X^* = F_X^{-1}(U)$, $Y^* = F_Y^{-1}(U)$ and $U \sim U(0,1)$.  
For $r > 1$, equality is attained only if $(X,Y) \sim (X^*,Y^*)$.
\end{lemma}
\begin{proof} Observe, as in 
\citeasnoun{frechet}, that if the random variables $X,Y$ have fixed
marginals $F_X$ and $F_Y$, then
\begin{equation} \label{eq:kupper}
\pr (X \leq x, Y \leq y) \leq H_+(x,y),\end{equation}
where $H_+(x,y) = \min(F_X(x),
F_Y(y))$. This bound is achieved by taking $U \sim 
U(0,1)$ and setting $X^* = F_X^{-1}(U), Y^* = F_Y^{-1}(U)$.

Thus, by Lemma \ref{lem:intprt}, with $k(x,y) = -|x-y|^r$, for $r \geq 1$,
and taking $H_1(x,y) = \pr(X \leq x, Y \leq y)$ and $H_2 = H_+$, we
deduce that
$$ \ep |X-Y|^r - \ep |X^*-Y^*|^r = \int \left\{ H_+(x,y) - H_1(x,y) \right\}
d\mu_k(x,y) \geq 0,$$
so $(X^*,Y^*)$ achieves the infimum in the definition of the Wasserstein
distance.

Finally, since taking $r >1$ implies that
the measure $\mu_k$ has a strictly positive density with respect to Lebesgue
measure,
we can only have equality in (\ref{eq:covmax}) if $\pr(X \leq x, Y \leq y) 
= \min \{ F_X(x), F_Y(y) \}$ Lebesgue almost everywhere.  
But the joint distribution function is right-continuous, 
so this condition determines the value of $\pr(X \leq x, Y \leq y)$ 
everywhere.   
\end{proof}
Using the construction in Lemma \ref{lem:covvalue}, \citeasnoun{bickel} 
establish that if
$X_1$ and $X_2$ are independent and $Y_1$ and $Y_2$ are independent, then
\begin{equation} \label{eq:sub}
d_2^2(X_1 + X_2, Y_1 + Y_2) \leq d_2^2(X_1, Y_1) + d_2^2(X_2, Y_2).
\end{equation}
Similar subadditive expressions arise  in the proof of convergence of  
Fisher information  in \citeasnoun{johnson5}. By  focusing on the
case $r=2$ in Definition \ref{def:mallows}, and by using the theory
of $L^2$ spaces and projections,  we establish parallels with the
Fisher information argument.

We prove Equation (\ref{eq:sub}) below, and further consider the
case of equality in this relation.
\citeasnoun[p.504]{major2} gives an equivalent construction to that
given in Lemma \ref{lem:covvalue}.
If  $F_Y$ is a continuous distribution function, then $F_Y(Y) \sim U(0,1)$,
so we generate $Y^* \sim F_Y$ and take $X^* = F_X^{-1} \circ F_Y(Y^*)$.
Recall that if $\ep X = \mu$ and $\var X = \sigma^2$, we
write $D(X)$ for $d_2^2(X,Z_{\mu,\sigma^2})$.
\begin{proposition} \label{prop:iff} If $X_1$, $X_2$ are independent,
with finite
variances $\sigma_1^2, \sigma_2^2 > 0$, then for any $t \in (0,1)$,
$$ D \left( \sqrt{t} X_1 + \sqrt{1-t} X_2 \right) \leq 
t D(X_1) + (1-t) D(X_2), $$
with equality if and only if 
$X_1$ and $X_2$ are normal.\end{proposition}
\begin{proof} We consider bounding $D(X_1 + X_2)$ for independent $X_1$
and $X_2$ with mean zero, since the general result follows on 
translation and rescaling.

We generate independent $Y_i^* \sim N(0,\sigma^2_i)$, and 
take $X_i^* = F_{X_i}^{-1}\circ \Phi_{\sigma^2_i}(Y_i^*) = h_i(Y_i^*)$, say,
for  $i=1,2$. 
Further, writing $\sigma^2 = \sigma_1^2 + \sigma^2_2$, 
we define $Y^* = Y^*_1 + Y_2^*$ and set
$X^* = F_{X_1 + X_2}^{-1} \circ \Phi_{\sigma^2}(Y^*_1 + Y^*_2) = 
h(Y_1^* + Y_2^*)$, say. Then
\begin{eqnarray*}
d_2^2( X_1 + X_2, Y_1 + Y_2) & = & \ep ( X^* - Y^*)^2 \\
& \leq & \ep (X_1^* + X_2^* - Y_1^* - Y_2^*)^2 \\
& = & \ep (X_1^* - Y_1^*)^2 + \ep (X_2^* - Y_2^*)^2 \\
& = & d_2^2(X_1, Y_1) + d_2^2(X_2, Y_2). \end{eqnarray*}
Equality holds if and only if 
$(X_1^* + X_2^*,Y_1^* + Y_2^*)$ has the same distribution as $(X^*,Y^*)$.
By our construction of $Y^* = Y_1^* + Y_2^*$, this means that 
$(X_1^*+ X_2^*, Y_1^* + Y_2^*)$ has the same distribution as 
$(X^*,Y^*_1+Y_2^*)$, so
$\pr \{ X_1^* + X_2^* = h(Y^*_1 + Y^*_2) \} = \pr \{ 
X^* = h(Y^*_1+Y_2^*) \} = 1$. 
Thus, if equality holds, then
\begin{equation} \label{eq:linear}
h_1(Y_1^*) + h_2(Y_2^*) = h(Y_1^* + Y_2^*) \mbox{ almost surely}.
\end{equation}
\citeasnoun{brown} and \citeasnoun{johnson5}, 
showed that equality holds in Equation (\ref{eq:linear}) 
if and only if $h, h_1, h_2$ are linear. In particular, Proposition 2.1 of
\cite{johnson5} implies that there exist 
constants $a_i$ and $b_i$ such that
\begin{eqnarray} \label{eq:brown}
\lefteqn{\ep \{ h(Y^*_1 + Y^*_2) - h_1(Y^*_1) - h_2(Y^*_2)\}^2} \nonumber \\
& \geq &  \frac{2\sigma_1^2 \sigma_2^2}{(\sigma_1^2 + \sigma_2^2)^2}
\left[ \ep \{h_1(Y^*_1) - a_1 Y^*_1 - b_1 \}^2 +
\ep \{ h_2(Y^*_2) - a_2 Y^*_2 - b_2 \}^2 \right]. \end{eqnarray}
Hence, if Equation (\ref{eq:linear})
holds, then $h_i(u) = a_i u + b_i$ almost everywhere. Since $Y_i^*$ and
$X_i^*$ have the same mean and variance, it follows that $a_i=1$, $b_i=0$.
Hence $h_1(u) = h_2(u) = u$ and $X_i^* = Y_i^*$.
\end{proof}
Recall that $T_k = S_{2^k}$, where $S_n = (X_1 + \ldots + X_n)/\sqrt{n}$
is a normalised sum of independent and identically distributed
 random variables of mean zero and finite
variance $\sigma^2$.

\begin{proposition} \label{prop:rich} 
There exists a strictly increasing sequence $(k_r) \in
\mathbb{N}$ and a random variable $T$ such that
\[
\lim_{r \rightarrow \infty} D(T_{k_r}) = D(T).
\]
If $T_{k_r}'$ and $T'$ are independent copies of $T_{k_r}$ and $T$
respectively, then
\[
\lim_{r \rightarrow \infty} D(T_{k_r+1}) =
\lim_{r \rightarrow \infty} D\biggl(\frac{T_{k_r} + T_{k_r}'}{\sqrt{2}}\biggr) 
= D\biggl(\frac{T + T'}{\sqrt{2}}\biggr).
\]
\end{proposition}
\begin{proof}
Since $\var(T_k) = 1$ for all $k$, the sequence $(T_k)$ is tight.
Therefore, by Prohorov's theorem, there exists a strictly increasing
sequence $(k_r)$ and a random variable $T$ such that
\begin{equation}
\label{weak}
T_{k_r} \stackrel{d}{\rightarrow} T
\end{equation}
as $r \rightarrow \infty$.  Moreover, the
proof of Lemma~5.2 of \citeasnoun{brown} shows that the sequence $(T_{k_r}^2)$
is uniformly integrable. 
But this, combined with Equation (\ref{weak})
implies that $\lim_{r \tends} d_2(T_{k_r},T) = 0$ 
\cite[Lemma 8.3(b)]{bickel}.  Hence
\[
D(T_{k_r}) = d_2^2(T_{k_r},Z_{\sigma^2}) \leq \{d_2(T_{k_r},T) +
d_2(T,Z_{\sigma^2})\}^2 \rightarrow d_2^2(T,Z_{\sigma^2}) = D(T)
\]
as $r \rightarrow \infty$.  Similarly, $d_2^2(T,Z_{\sigma^2}) \leq
\{d_2(T,T_{k_r}) + d_2(T_{k_r},Z_{\sigma^2})\}^2$, yielding the
opposite inequality.  This proves the first part of the proposition.

For the second part, it suffices to observe that $T_{k_r} + T_{k_r}'
\stackrel{d}{\rightarrow} T + T'$ as $r \rightarrow \infty$, and
$\ep (T_{k_r} + T_{k_r}')^2 \rightarrow \ep(T + T')^2 $,
and then use the same argument as in the first part of the proposition.
\end{proof}
Combining Propositions \ref{prop:iff} and \ref{prop:rich}, as described
in Section \ref{sec:intro}, the proof of Theorem \ref{thm:main} is now
complete.
\section{Convergence of $d_r$ for general $r$}
The subadditive inequality (\ref{eq:sub}) arises in part
from a moment inequality; that is, if
$X_1$ and $X_2$ are independent with mean zero,
then $\ep |X_1 + X_2|^r \leq \ep |X_1|^r + \ep |X_2|^r$, for $r=2$. Similar
results imply that for $r \geq 2$, we have
$\lim_{n \tends} d_r(S_n, Z_{\sigma^2}) = 0$. 
First, we prove the following lemma:
\begin{lemma} \label{lem:rsub}
Consider independent random variables $V_1, V_2, \ldots$ and 
$W_1, W_2, \ldots$, where for some $r \geq 2$ and for all $i$,
$\ep |V_i|^r < \infty$ and $\ep |W_i|^r < \infty$.
Then for any $m$, there exists
a constant $c(r)$ such that
\begin{eqnarray*}
\lefteqn{d_r^r \left( V_1 + \ldots + V_m, W_1 + \ldots + W_m \right) }\\
& \leq & c(r) \biggl\{ \sum_{i=1}^m d_r^r(V_i, W_i) + 
\biggl( \sum_{i=1}^m d_2^2(V_i, W_i) \biggr)^{r/2} \biggr\}.
\end{eqnarray*} 
\end{lemma}
\begin{proof}
We consider independent $U_i \sim U(0,1)$, and set
$V_i^* = F_{V}^{-1}(U_i)$ and $W_i^* = F_{W}^{-1}(U_i)$. Then
\begin{eqnarray*}
\lefteqn{d_r^r(V_1 + \ldots + V_m, W_1 + \ldots + W_m) } \\
& \leq & \ep \left| \sum_{i=1}^m (V_i^* - W_i^*) \right|^r \nonumber \\
& \leq & c(r) \biggl\{ \sum_{i=1}^m \ep \left| V_i^* - W_i^* \right|^r 
+ \biggl( \sum_{i=1}^m \ep \left| V_i^* - W_i^* \right|^2 \biggr)^{r/2} 
\biggr\} \end{eqnarray*}
as required. This final line is an application of Rosenthal's inequality 
\cite[Theorem 2.9]{petrov} to the sequence $(V_i^* - W_i^*)$.
\end{proof}
Using Lemma \ref{lem:rsub}, we establish the following theorem. 
\begin{theorem} \label{thm:mainr}
Let $X_1, X_2, \ldots $ be 
independent and identically distributed random variables
with mean zero, variance $\sigma^2 > 0$ and 
$\ep |X_1|^r < \infty$ for some $r \geq 2$. If
$S_n = (X_1 + \ldots + X_n)/\sqrt{n}$, then
$$ \lim_{n \tends} d_r(S_n, Z_{\sigma^2}) = 0.$$
\end{theorem}
\begin{proof}
Theorem \ref{thm:main} covers
the case of $r=2$, so need only consider $r >2$.
We use a scaled version of 
Lemma \ref{lem:rsub} twice. First, we use $V_i = X_i$, 
$W_i \sim N(0, \sigma^2)$ and $m = n$, in order to deduce that, by monotonicity
of the $r$-norms:
\begin{eqnarray*}
 d_r^r \left( S_n, Z_{\sigma^2} \right) 
& \leq & 
c(r) \left\{ n^{1 - r/2} d_r^r(X_1, Z_{\sigma^2}) 
+ d_2^2(X_1, Z_{\sigma^2})^{r/2} \right\} \\
& \leq & c(r) \left( n^{1 - r/2} +1 \right) d_r^r(X_1, Z_{\sigma^2}),
\end{eqnarray*}
so that $d_r^r \left( S_n, Z_{\sigma^2} \right)$ is uniformly bounded in $n$,
by $K$, say.
Then, for general $n$, define $N = \lceil \sqrt{n} \rceil$, take
$m = \lceil n/N \rceil$, 
and $u = n - (m-1)N \leq N$. In Lemma \ref{lem:rsub}, take
\begin{eqnarray*}
V_i & = & X_{(i-1)N + 1} + \ldots + X_{i N}, 
\mbox{ for $i = 1, \ldots, m-1$} \\
V_m & = & X_{(m-1)N + 1} + \ldots + X_{n},  
\end{eqnarray*}
and $W_i \sim N(0, N \sigma^2)$ for $i = 1, \ldots, m-1$, 
$W_m \sim N(0, u \sigma^2)$ independently.
Now the uniform bound above gives, on rescaling,
$$d_{r}^r(V_i,W_i) = N^{r/2} d_r^r(S_N, Z_{\sigma^2})
\leq N^{r/2} K \mbox{ for $i = 1, \ldots m-1$} $$ and 
$d_r^r(V_m,W_m) = u^{r/2} d_r^r(S_u, Z_{\sigma^2}) \leq N^{r/2} K$.
Further $d_{2}^2(V_i,W_i) = N d_2^2(S_N, Z_{\sigma^2})$
 for $i = 1, \ldots m-1$ and 
$d_2^2(V_m,W_m) = u d_2^2(S_u, Z_{\sigma^2}) \leq N d_2^2(S_1, Z_{\sigma^2})$.
Hence, using Lemma \ref{lem:rsub} again, we obtain
\begin{eqnarray*}
\lefteqn{d_r^r \left( S_{n}, Z_{\sigma^2} \right) } \nonumber \\
& = & \frac{1}{n^{r/2}} 
d_r^r \left( V_1 + \ldots + V_m, W_1 + \ldots + W_m \right) \nonumber \\
& \leq & 
\frac{c(r)}{n^{r/2}} \left\{ \sum_{i=1}^m d_r^r(V_i, W_i) +
\left( \sum_{i=1}^m d_2^2(V_i, W_i) \right)^{r/2} \right\} \nonumber \\
& \leq &
c(r) \left\{ m K \frac{N^{r/2}}{n^{r/2}} +
\left( \frac{N(m-1)}{n} d_2^2(S_N, Z_{\sigma^2}) + \frac{N}{n} 
d_2^2(S_1, Z_{\sigma^2}) \right)^{r/2} \right\} \nonumber \\
& \leq &
c(r) \left\{ \frac{m K}{(m-1)^{r/2}} +
\left( d_2^2(S_N, Z_{\sigma^2}) + \frac{1}{m-1} 
d_2^2(S_1, Z_{\sigma^2}) \right)^{r/2} \right\}. 
\label{eq:toconv}
\end{eqnarray*}
This
converges to zero since $\lim_{n \tends} d_2(S_N, Z_{\sigma^2}) = 0$.
\end{proof}
\section{Strengthening subadditivity}
Under certain conditions, we obtain a rate for the convergence in
Theorem \ref{thm:main}. 
Equation (\ref{eq:subadd})  shows that
$D(T_{k})$  is decreasing.  Since  $D(T_{k})$ is  bounded
below, the difference  sequence $D(T_k) - D(T_{k+1})$ 
converges  to zero,  As in  \citeasnoun{johnson5}  we examine
this difference sequence, to show  that its convergence
implies  convergence of $D(T_{k})$ to zero.

Further, in the  spirit of \citeasnoun{johnson5}, 
we hope  that if the difference
sequence  is small,  then  equality  `nearly'  holds in Equation
(\ref{eq:linear}),  and so  the
functions $h,  h_1, h_2$ are `nearly'  linear. 
This implies that if $\cov(X,Y)$ is
close to its  maximum, then $X$ is be close to  $h(Y)$ in the $L^2$ sense.

Following \citeasnoun{delbarrio},
we define a new distance quantity 
$ D^*(X) = \inf_{m,s^2}
d_2^2(X, Z_{m,s^2}). $
Notice that $D(X) = 2 \sigma^2 - 2 \sigma k \leq 2 \sigma^2$, where
$k = \int_0^1 F_X^{-1}(x) \Phi^{-1}(x) dx$.
This follows since $F_X^{-1}$ and $\Phi^{-1}$ are increasing functions, so
$k \geq 0$ by Chebyshev's rearrangement lemma. 
Using results of \citeasnoun{delbarrio}, it follows that
$$ D^*(X) = \sigma^2 - k^2 = D(X) - \frac{D(X)^2}{4 \sigma^2},$$
and convergence of $D(S_n)$ to zero is 
equivalent to convergence of $D^*(S_n)$ to zero.
\begin{proposition} \label{prop:double}
Let $X_1$ and $X_2$ be independent and identically distributed
random variables with mean $\mu$, variance $\sigma^2 > 0$ and
densities (with respect to Lebesgue measure).
Defining $g(u) =  \Phi_{\mu,\sigma^2}^{-1} \circ 
F_{(X_1 + X_2)/\sqrt{2}}(u)$, if the
derivative $g'(u) \geq c$ for all $u$ then
$$ D \left( \frac{ X_1 + X_2}{\sqrt{2}} \right) 
\leq \left( 1 - \frac{c}{2} \right) D(X_1) + \frac{c D(X_1)^2}{8 \sigma^2}
\leq
\left( 1 - \frac{c}{4} \right) D(X_1).$$
\end{proposition}
\begin{proof} As before, translation invariance allows us to take $\ep X_i 
=0$.
For random variables $X,Y$,
we    consider      the     difference     term     Equation
(\ref{eq:kupper}) and write $g(u) = F_Y^{-1} \circ F_X(u)$, and $h(u)
= g^{-1}(u)$. The function $k(x,y) = -\{ x - h(y) \}^2$ is 
quasi-monotone and induces the measure $d\mu_k(x,y) = 2 h'(y) dx dy$. Taking 
$H_1(x,y) = \pr(X \leq x, Y \leq y)$ and $H_2(x,y) = \min \{ F_X(x),
F_Y(y) \}$ in  
Lemma \ref{lem:intprt} implies that
$$ \ep \{X - h(Y) \}^2 = 2 \int h'(y) \left\{ H_2(x,y) - H_1(x,y) \right\}
dx dy,$$
since $\ep \{X^* - h(Y^*) \}^2 = 0$. By assumption $h'(y) \leq 1/c$,
so 
$$ \ep \{ X - h(Y) \}^2 \leq \frac{2}{c} \left\{ \cov(X^*, Y^*) - 
\cov(X,Y) \right)\}.$$  
Again take $Y_1^*,Y_2^*$ independent $N(0,\sigma^2)$ and set
$X_i^* = F_{X_i}^{-1} \circ F_{Y_i}(Y_i^*) = h_i(Y_i^*)$. Then define
$Y^* = Y_1^* + Y_2^*$ and take $X^* = F_{X_1 + X_2}^{-1} \circ F_{Y_1 + Y_2}
(Y^*)$. Then there exist $a$ and $b$ such that
\begin{eqnarray*}
\lefteqn{d_2^2(X_1, Y_1) + d_2^2(X_2, Y_2) - d_2^2(X_1 + X_2, Y_1 + Y_2)} \\
& = & \ep(X_1^* + X_2^* - Y_1^* - Y_2^*)^2 - \ep(X^* - Y^*)^2 \\
& = & 2 \cov(X^*,Y^*) - 2 \cov(X_1^* + X_2^*,Y_1^* + Y_2^*) \\
& \geq & c \ep \{ X_1^* + X_2^* - h(Y^*_1 + Y^*_2) \}^2 \\
& = & c \ep \{h_1(Y_1^*) + h_2(Y_2^*) - h(Y_1^* + Y_2^*)\}^2 \\
& \geq & c \ep \{ h_1(Y_1^*) - a Y_1^* - b \}^2 \geq c D^*(X_1), 
\end{eqnarray*}
where the penultimate inequality follows by Equation (\ref{eq:brown}).
Recall that $D(X) \leq 2 \sigma^2$, so that 
$D^*(X) = D(X) - D(X)^2/(4 \sigma^2)
\geq D(X)/2$. The result follows on rescaling.
\end{proof}
We briefly discuss the strength of the condition imposed. If
 $X$ has mean zero, distribution function $F_X$ and continuous
density $f_X$, define the scale invariant quantity
$$ \C(X) = \inf_{u} (\Phi_{\sigma^2}^{-1} \circ F_X)'(u)
= \inf_{p \in (0,1)} \frac{ f_X(F_{X}^{-1}(p))}
{ \phi_{\sigma^2}( \Phi_{\sigma^2}^{-1}(p))}
= \inf_{p \in (0,1)} \sigma \frac{ f_X(F_{X}^{-1}(p))}
{ \phi( \Phi^{-1}(p))}.$$
We want to understand when $\C(X) > 0$. 
\begin{example} If $X \sim U(0,1)$, then
$\C(X) = 1/\sqrt{12 \sup_x \phi(x)} = \sqrt{\pi/6}.$ 
\end{example}
\begin{lemma} If $X$ has mean zero and variance $\sigma^2$ then
$\C(X)^2 \leq \sigma^2/(\sigma^2 + {\rm median}(X)^2)$. \end{lemma}
\begin{proof}
By the Mean Value Inequality, for all $p$ 
$$ |\Phi_{\sigma^2}^{-1}(p)| = | \Phi_{\sigma^2}^{-1}(p) - 
\Phi_{\sigma^2}^{-1}(1/2) | 
\geq \C(X) | F_X^{-1}(p) - F_X^{-1}(1/2) |,$$
so that 
\begin{align*}
\sigma^2 + F^{-1}_X(1/2)^2 & =  \int_0^1 F_X^{-1}(p)^2 dp + 
F^{-1}_X(1/2)^2  
=  \int_0^1 \{ F_X^{-1}(p) - F^{-1}_X(1/2) \}^2 dp  \\
& \leq  \frac{1}{\C(X)^2} \int_0^1 \Phi_{\sigma^2}^{-1}(p)^2 dp  
=  \frac{\sigma^2}{\C(X)^2}.
\end{align*}
\end{proof}
In general we are concerned with the rate at which $f_X(x) \rightarrow 
0$ at the edges of the support.
\begin{lemma}
If for some $\epsilon > 0$, 
\begin{equation} \label{eq:denspow2}
f_X(F_X^{-1}(p)) \simeq c (1-p)^{1 - \epsilon} 
\mbox{ as $p \rightarrow 1$}\end{equation} then
$\lim_{p \rightarrow 1} f_X(F_X^{-1}(p))/\phi(\Phi^{-1}(p)) = 
\infty$.
Correspondingly if 
\begin{equation} \label{eq:denspow} f_X(F_X^{-1}(p)) \simeq c p^{1 - \epsilon} 
\mbox{ as $p \rightarrow 0$} \end{equation} then
$ \lim_{p \rightarrow 0} f_X(F_X^{-1}(p))/\phi(\Phi^{-1}(p)) = 
\infty.$
\end{lemma}
\begin{proof} 
Simply note that by the Mills ratio 
\cite[p.850]{shorack} as $x \rightarrow \infty$, 
$\Phi(x) \sim \phi(x)/ x$, so that as $p \rightarrow 1$,
$\phi(\Phi^{-1}(p)) \sim (1-p) \Phi^{-1}(p) 
\sim (1-p) \sqrt{-2 \log (1-p)}$.
\end{proof}
\begin{example} \mbox{   } \; \; 

\begin{enumerate}
\item{The density of the $n$-fold convolution of $U(0,1)$
random variables is given by $f_X(x) = 
x^{n-1}/(n-1)!$ for $0 < x < 1$, hence $F_X^{-1}(p) 
= (n! p)^{1/n}$, and
$f_X(F_X^{-1}(p)) = n/(n!)^{1/n} p^{(n-1)/n}$, so that 
Equation (\ref{eq:denspow}) holds.}
\item{For an ${\rm Exp}(1)$ random variable, $f_X(F_X^{-1}(p)) = 1- p$, so that
Equation (\ref{eq:denspow2}) fails and $\C(X) = 0$.}
\end{enumerate}
\end{example}
To obtain bounds on $D(S_n)$ as $n \tends$, we need to control
the  sequence $\C(S_n)$.  Motivated  by properties  of the  (seemingly
related)  Poincar\'{e}  constant,  we  conjecture  that  $\C(  (X_1  +
X_2)/\sqrt{2}) \geq \C(X_1)$ for independent and identically distributed $X_i$.
If this is true and $\C(X) = c$ then $\C(S_n) \geq c$ for all $n$.

Assuming that $\C(S_n) \geq c$ for all $n$, note that $D(T_k) \leq 
(1 - c/4)^k D(X_1) \leq (1 - c/4)^k (2 \sigma^2)$.  Now
$$ D(T_{k+1}) \leq D(T_k) (1 -c/2) \left\{ 1 + \frac{ c D(T_k)}{8\sigma^2
(1-c/2)} \right\},$$
so 
$$\prod_{k=0}^{\infty} 
\left\{ 1 + \frac{c D(T_k)}{ 8\sigma^2 (1-c/2)} \right\}
 \leq \exp \left\{ \sum_{k=0}^{\infty} \frac{c D(T_k)}{8 \sigma^2(1-c/2)} 
\right\} 
\leq \exp \left( \frac{1}{1-c/2} \right).$$
We deduce that 
$$ D(T_{k})  \leq D(X_1) \exp \left( \frac{1}{1-c/2} \right) 
\left( 1- c/2 \right)^k ,$$ or  
$D(S_n) = O(n^t)$, where $t = \log_2(1-c/2).$
\begin{remark} In general, convergence of $d_4(S_n, Z_{\sigma^2})$ cannot
occur at a rate faster than $O(1/n)$. This follows because  
$\ep S_n^4 = 3 \sigma^4 + \gamma(X_1)/n$,
where $\gamma(X)$, the excess kurtosis, is defined by $\gamma(X) = \ep X^4
- 3( \ep X^2)^2$ (when $\ep X =0$). Thus by Minkowski's inequality,
\begin{eqnarray*}
d_4(S_n,Z_{\sigma^2}) & \geq & \left| (\ep S_n^4)^{1/4} - 
(\ep Z_{\sigma^2}^4)^{1/4} \right| \\
& = & 3^{1/4} \sigma \left| 
\left(1 + \frac{\gamma(X)}{n} \right)^{1/4} - 1 \right| 
= \frac{3^{1/4} \sigma |\gamma(X)|}{4 n} + O\left( \frac{1}{n^2} \right).
\end{eqnarray*}
\end{remark}
Motivated  by  this remark, and by analogy with the rates discovered in 
\citeasnoun{johnson5},
we conjecture  that the true rate of convergence
is $D(S_n)  = O(1/n)$.  To obtain  this, we would need to control 
$1-\C(S_n)$. 
\section{Convergence to stable distributions} \label{sec:stable}
We now consider convergence to other stable 
distributions. \citeasnoun{gnedenko} review
classical results of this kind. We say that $Y$ is $\alpha$-stable
if, when $Y_1, \ldots Y_n$ are independent copies of $Y$, we have
$ (Y_1 + \ldots + Y_n - b_n)/n^{1/\alpha} \sim Y$ for some sequence $(b_n)$.
Note that $\alpha$-stable variables only exist for $0 < \alpha
\leq 2$; we assume for the rest of this Section that $\alpha < 2$. 
\begin{definition} \label{def:dna}
If $X$ has a distribution function of the form
\begin{eqnarray*}
F_X(x) & = & \frac{c_1 + b_X(x)}{|x|^{\alpha}} \mbox{ for $x < 0$} \\
1 - F_X(x) & = & \frac{c_2 + b_X(x)}{x^{\alpha}} \mbox{ for $x \geq 0$}
\end{eqnarray*}
where $b_X(x) \rightarrow 0$ as $x \rightarrow \pm \infty$, then we say that
$X$ is in the domain of normal attraction of some stable $Y$ with tail
parameters $c_1,c_2$. 
\end{definition}
Theorem 5 of Section 35 of \cite{gnedenko} shows that if $F_X$ is
of this form, there exist a sequence $(a_n)$ and
an $\alpha$-stable distribution function $F_Y$, 
determined by the parameters $\alpha$, $c_1$, $c_2$, such that 
\begin{equation} \label{eq:stabconv} \frac{ X_1 + \ldots + X_n - 
a_n}{n^{1/\alpha}} 
\stackrel{d}{\rightarrow} F_Y.\end{equation}

Although Equation (\ref{eq:stabconv}) 
is obviously very similar to the standard Central Limit Theorem,
one important distinguishing feature is that both $\ep |X|^{\alpha}$
and $\ep |Y|^{\alpha}$ are infinite for $0 < \alpha < 2$.

We use the following moment bounds from
\citeasnoun{vonbahr2}. If $X_1, X_2, \ldots$ are independent, then
\begin{eqnarray}
\ep | X_1 + \ldots + X_n |^r & \leq & \sum_{i=1}^n \ep |X_i|^r
\quad \mbox{for $0 < r \leq 1$} \label{eq:vb1} \\
\ep | X_1 + \ldots + X_n |^r & \leq & 2 \sum_{i=1}^n \ep |X_i|^r
\quad \parbox{1.2in}{when $\ep X_i =0$, 
for $1 < r \leq 2$.}
\label{eq:vb2}
\end{eqnarray}
Now, using ideas of \citeasnoun{stout}, we show that for a subset of the
domain of normal attraction, $d_{\beta}(X,Y) < \infty$, for some
$\beta > \alpha$.
\begin{definition} \label{def:dsna}
We say that a random variable is in the 
domain of strong normal attraction of $Y$ 
if the function $b_X(x)$ from Definition \ref{def:dna} satisfies
$$ b_X(x) \leq \frac{C}{|x|^{\gamma}},$$
for some constant $C$ and some $\gamma > 0$. \end{definition}
\citeasnoun{cramer} shows that such random variables have an
Edgeworth-style expansion, and thus 
convergence to $Y$ occurs. However, his proof requires some
involved analysis and use of characteristic functions. See also 
\citeasnoun{mijnheer} and \citeasnoun{mijnheer2}, which use
bounds based on the quantile transformation described above.

We can regard Definition \ref{def:dsna} as being 
analogous to requiring 
a bounded $(2+\delta)$th moment in the Central Limit Theorem,
which allows an explicit rate of convergence (via
the Berry-Ess\'{e}en theorem).
We now show the relevance of Definition \ref{def:dsna} to the problem of stable
convergence.
\begin{lemma} \label{lem:dsna}
If $X$ is in the domain of strong normal attraction of
an $\alpha$-stable random variable $Y$,
then $d_{\beta}(X,Y) < \infty$ for some $\beta > \alpha$. \end{lemma}
\begin{proof} 
We show that Major's construction always gives a joint distribution 
$(X^*,W^*)$ with $\ep |X^* - W^*|^{\beta} < \infty$, and hence
$d_{\beta}(X,W) < \infty$. 
Following \citeasnoun{stout}, define a random variable
$W$ by
\begin{eqnarray*} 
\pr(W \geq x) & = & 
c_2 x^{-\alpha} \mbox{  if $x > (2c_2)^{1/\alpha}$.} \\
\pr(W \leq x) & = &  c_1 |x|^{-\alpha}  \mbox{  
if $x < -(2c_1)^{1/\alpha}$.} \\
\pr(W \in [-(2c_1)^{1/\alpha}, (2c_2)^{1/\alpha}]) & = & 0.
\end{eqnarray*} 
Then for $w > 1/2$, $F_W^{-1}(w) = \{ c_2/(1-w) \}^{1/\alpha}$, and so for
$x \geq 0$,
$$ x - F_W^{-1} (F_X(x)) = x \left\{ 1 - \left( \frac{c_2}{c_2 + b_X(x)} 
\right)^{1/\alpha} \right\}.$$
Now, since $b_X(x) \rightarrow 0$, there exists $K$ such that 
if $x \geq K$ then $b_X(x) \geq - c_2/2$.

By the Mean Value Inequality, if $t \geq -1/2$, then
$$ \left| 1 - (1 + t)^{-1/\alpha} \right| \leq \frac{ |t| 2^{1+1/\alpha}}{
\alpha},$$
so that for $x \geq K$
$$ \left| x - F_W^{-1} F_X(x) \right| 
\leq \frac{2^{1 + 1/\alpha} x |b_X(x)|}{\alpha c_2}.$$
Thus, if $X$ is in the strong domain of attraction, then
$$ \int_{|x| \geq K} \left| x - F_W^{-1} F_X(x) \right|^{\beta} 
dF_X(x) \leq \left(  \frac{2^{1 + 1/\alpha} C}{\alpha c_2} \right)^{\beta}
\int_{|x| \geq K} |x|^{\beta(1- \gamma)}  dF_X(x).$$
Hence $d_{\beta}(X,W)$ is finite for all $\beta$ if 
$\gamma \geq 1$ and for $\beta < \alpha/(1- \gamma)$, if $\gamma < 1$.
Moreover,
\citeasnoun[Equation (2.2)]{mijnheer2} shows that if $Y$ is $\alpha$-stable,
then as $x \tends$,
$$ \pr(Y \geq x) = \frac{c_2}{x^{\alpha}} + 
O \left( \frac{1}{x^{2\alpha}} \right).$$
and so $Y$ is in its own domain of strong normal attraction. 
Thus using the construction above, $d_{\beta}(Y,W)$ is finite for all
$\beta$ if $\alpha \geq 1$ and for $\beta < \alpha/(1-\alpha)$ otherwise.

Recall that the triangle inequality holds, for $d_{\beta}$ or
$d_{\beta}^{\beta}$, according as $\beta \geq 1$ or $\beta < 1$.
Hence $d_{\beta}(X,Y)$ is finite for all $\beta$ if
$\min(\alpha,\gamma) \geq 1$ and for 
$\beta < \alpha/(1 - \min(\alpha,\gamma))$ otherwise.
\end{proof}
Note that for random variables $X_i$ in the same strong domain of normal
attraction,
$d_{\beta}(X_i,Y)$ may be bounded in terms of the function $b_{X_i}(x)$.
In particular
if there exist $C,\gamma$ such that $b_{X_i}(x) \leq C/|x|^{\gamma}$
then $\sup_i d_{\beta}(X_i,Y) < \infty$,
so the hypothesis of Theorem \ref{thm:mainstab} is satisfied.

\begin{proof}{\bf of Theorem \ref{thm:mainstab}}
We use the bounds provided by Equations (\ref{eq:vb1}) and
(\ref{eq:vb2}). We consider independent pairs $(X_i^*,Y_i^*)$ having the joint
distribution that achieves the infimum in Definition \ref{def:mallows}.
Then by rescaling we have that
\begin{eqnarray*}
d_{\beta}^{\beta}(S_n, Y) 
& \leq & \frac{1}{n^{\beta/\alpha}} d_{\beta}^{\beta}(X_1 + \ldots 
+ X_n, Y_1 + \ldots + Y_n)  \\
& \leq & \frac{1}{n^{\beta/\alpha}}
\ep \biggl| \sum_{i=1}^n (X_i^* - Y_i^*) \biggr|^{\beta} 
\leq  \frac{2}{n^{\beta/\alpha}}
\sum_{i=1}^n  \ep \left| X_i^* - Y_i^* \right|^{\beta}. 
\end{eqnarray*}
We deduce that in the case of identical variables, 
$d_{\beta}(S_n,Y)$ (and hence 
$d_{\alpha}(S_n,Y)$) converges at rate $O(n^{1/\beta - 1/\alpha})$.
\end{proof}
We now combine Theorem \ref{thm:mainstab} and Lemma \ref{lem:dsna}, to
obtain a rate of convergence for identical variables. 
Note that Theorem \ref{thm:mainstab} 
requires us to take $\beta \leq 2$. Overall then
we deduce that $d_{\alpha}(S_n,Y)$ converges at rate $O(n^{-t})$, where
\begin{enumerate}
\item{if $\min(\alpha, \gamma) \geq 1$, we take $\beta =2$, and hence
$t = 1/\alpha - 1/2$;} 
\item{if $\min(\alpha,\gamma) < 1$, we may take 
$\beta = \min[ \alpha/\{1 - \min(\alpha,\gamma) + \epsilon\}, 2]$
for any $\epsilon >0$, and then 
$t = \min(1/\alpha - 1/2, 1 - \epsilon, \gamma/\alpha-\epsilon) $.}
\end{enumerate}

Theorem \ref{thm:mainr} implies that
if $d_r(S_n, Z_{\sigma^2})$ ever becomes finite, then it tends to zero,
the counterpart of the following result.
\begin{theorem} \label{thm:stablex}
Fix $\alpha \in (0,2)$, let $X_1, X_2, \ldots $ be 
independent random 
variables (where $\ep X_i =0$, if $\alpha > 1$), and let 
$S_n = (X_1 + \ldots + X_n)/n^{1/\alpha}$.  
Suppose there exists an $\alpha$-stable random variable $Y$ and
$Y_1, Y_2, \ldots$ having the same distribution as $Y$, and satisfying
\begin{equation} \label{eq:linde}
\frac{1}{n} \sum_{i=1}^n \ep \{ |X_i - Y_i|^{\alpha} \I( |X_i - Y_i| 
> b) \} \rightarrow 0 \mbox{ as $b \tends$.}\end{equation}
If $\alpha \neq 1$ then $\lim_{n \tends} d_{\alpha}(S_n, Y) = 0$,
and if $\alpha = 1$ then there exists a sequence
$c_n = n^{-1} \sum_{i=1}^n \ep(X_i - Y_i)$ such that
$\lim_{n \tends} d_{\alpha}(S_n -c_n, Y) = 0$.
\end{theorem}
\begin{proof}(Suggested by an anonymous referee).
Fix $\epsilon  > 0$. Suppose first that $1 \leq \alpha < 2$ and let
$d_i = \ep(X_i - Y_i)$. Note that $d_i=0$ if $\alpha > 1$. 
Let $b > 0$ and define
\begin{eqnarray*}
U_i & = & (X_i - Y_i) \I( |X_i - Y_i| \leq b) -
\ep \{ (X_i - Y_i) \I( |X_i - Y_i| \leq b) \} \\
V_i & = & (X_i - Y_i) \I( |X_i - Y_i| > b) -
\ep \{ (X_i - Y_i) \I( |X_i - Y_i| > b) \}.
\end{eqnarray*}
Then by Equation (\ref{eq:vb2}),
\begin{align*}
d^{\alpha}_{\alpha}(S_n - c_n,Y)
& \leq  \frac{1}{n} \ep \biggl| \sum_{i=1}^n (X_i - Y_i - d_i) 
\biggr|^{\alpha}
=  \frac{1}{n} \ep \biggl| \sum_{i=1}^n U_i + \sum_{i=1}^n V_i 
\biggr|^{\alpha} \\
& \leq  \frac{2^{\alpha-1}}{n} \ep \biggl| \sum_{i=1}^n U_i \biggr|^{\alpha} 
+ \frac{2^{\alpha-1}}{n} \ep \biggl| \sum_{i=1}^n V_i \biggr|^{\alpha} \\
& \leq  \frac{2^{\alpha-1}}{n} \biggl\{ \ep \biggl( \sum_{i=1}^n U_i \biggr)^2
\biggr\}^{\alpha/2} + 
\frac{2^{\alpha}}{n} \sum_{i=1}^n \ep | V_i |^{\alpha} \\
& \leq  \frac{2^{\alpha-1}}{n} 
\biggl( \sum_{i=1}^n \ep U_i^2 \biggr)^{\alpha/2} +
 \frac{2^{2 \alpha -1}}{n} \sum_{i=1}^n 
\ep \{ |X_i - Y_i|^{\alpha} \I( |X_i - Y_i| > b) \} \\
& \qquad \qquad + \frac{2^{2\alpha-1}}{n} \sum_{i=1}^n
\left[ \ep \left\{ |X_i - Y_i| \I( |X_i - Y_i| > b) \right\} 
\right]^{\alpha} \\ 
& \leq  \frac{2^{\alpha-1} b^{\alpha}}{n^{1-\alpha/2}} 
+ \frac{2^{2 \alpha}}{n} \sum_{i=1}^n 
\ep \{ |X_i - Y_i|^{\alpha} \I( |X_i - Y_i| > b) \} 
\end{align*}
The result follows on
choosing $b$ sufficiently large to control the second term, and then
$n$ sufficiently large to control the first.

For $0 < \alpha < 1$, take $U_i$ as before, take
$V_i =  (X_i - Y_i) \I( |X_i - Y_i| > b)$ and
$a_i = \ep \{ (X_i - Y_i) \I( |X_i - Y_i| \leq b) \}$. Now using 
Equation (\ref{eq:vb1}),
\begin{eqnarray*}
d^{\alpha}_{\alpha}(S_n, Y) 
& \leq & 
\frac{1}{n} \ep \biggl| \sum_{i=1}^n U_i + \sum_{i=1}^n V_i + \sum_{i=1}^n a_i 
\biggr|^{\alpha} \\
& \leq & 
\frac{1}{n} \ep \biggl| \sum_{i=1}^n U_i \biggr|^{\alpha} 
+ \frac{1}{n} \ep \biggl| \sum_{i=1}^n V_i \biggr|^{\alpha} 
+ \frac{1}{n} \biggl| \sum_{i=1}^n a_i \biggr|^{\alpha}  \\
& \leq & 
\frac{1}{n} \biggl\{ \ep \biggl( \sum_{i=1}^n U_i \biggr)^2 
\biggr\}^{\alpha/2}  
+ \frac{1}{n} \sum_{i=1}^n \ep |V_i|^{\alpha} 
+ \frac{b^{\alpha}}{n^{1-\alpha}},
\end{eqnarray*}
so again since $b$ is arbitrary, the result follows.
\end{proof}
Note when $X_1, X_2, \ldots$ are identically distributed, the Lindeberg
condition (\ref{eq:linde}) reduces to the requirement that $d_{\alpha}(X_1,Y)
< \infty$.

\setlength{\parskip}{0pt}

\end{document}